\documentclass[pos]{birkmult}
\usepackage{amsmath, amsthm,
amssymb}
 \newtheorem{thm}{Theorem}

\renewcommand{\le}{\leqslant}
\renewcommand{\ge}{\geqslant}

\newcommand{\diag}%
{\mathop{\rm
diag}\nolimits}

\newcommand{\rank}%
{\mathop{\rm
rank}\nolimits}

\newcommand{\ci}{
\begin{picture}(4,4)
\put(2,2){\circle*{2}}
\end{picture}}

\begin{document}
%
%
%
%
%
%
%
%
%
\title[Positivity
criteria generalizing
the leading principal
minors
criterion]{Positivity
Criteria Generalizing
the Leading Principal
Minors Criterion}


\author[Vyacheslav Futorny]{%
Vyacheslav Futorny}

\address{%
Department of
Mathematics,
University of S\~{a}o
Paulo,\\ Caixa Postal
6681, S\~{a}o Paulo,
Brazil}

\email{futorny@ime.usp.br}

\thanks{This is the authors' version of a work that was published in \emph{Positivity}, 11 (no. 1) (2007) 191--199.
The first author
was partially
supported by CNPq,
processo
307812/2004-9. The
second author was
partially supported by
FAPESP (S\~ao Paulo),
processo 05/59407-6.}


\author[Vladimir
Sergeichuk]{Vladimir
V. Sergeichuk}

\address{%
Institute of
Mathematics,\\
Tereshchenkivska 3,
Kiev, Ukraine}

\email{sergeich@imath.kiev.ua}


\author[Nadya
Zharko]{Nadya Zharko}
\address{Mech.-Math.
Faculty,
Kiev National University,\\
Vladimirskaya 64,
Kiev, Ukraine}

\email{n.zharko@mail.ru}



\date{January 1, 2004}


\begin{abstract}
An $n\times n$
Hermitian matrix is
positive definite if
and only if all
leading principal
minors
$\Delta_1,\dots,\Delta_n$
are positive. We show
that certain sums
$\delta_l$ of $l\times
l$ principal minors
can be used instead of
$\Delta_l$ in this
criterion. We describe
all suitable sums
$\delta_l$ for
$3\times 3$ Hermitian
matrices. For an
$n\times n$ Hermitian
matrix $A$ partitioned
into blocks $A_{ij}$
with square diagonal
blocks, we prove that
$A$ is positive
definite if and only
if the following
numbers $\sigma_l$ are
positive: $\sigma_l$
is the sum of all
$l\times l$ principal
minors that contain
the leading block
submatrix
$[A_{ij}]_{i,j=1}^{k-1}$
(if $k>1$) and that
are contained in
$[A_{ij}]_{i,j=1}^{k}$,
where $k$ is the index
of the block $A_{kk}$
containing the $(l,l)$
diagonal entry of $A$.
We also prove that
$\sigma_l$ can be used
instead of $\Delta_l$
in other inertia
problems.
\end{abstract}
\maketitle


\section{Introduction}

Let $A=[a_{ij}]$ be an
$n\times n$ Hermitian
complex matrix. By the
leading principal
minors criterion, $A$
is positive definite
if and only if
\begin{equation}\label{iytd1}
\Delta_1:=a_{11}>0,\quad
\Delta_2:=\begin{vmatrix}
a_{11}&a_{12}\\a_{21}&a_{22}
\end{vmatrix}>0,\ \dots,\
\Delta_n:=\begin{vmatrix}
a_{11}&\dots&a_{1n}\\
\vdots&\ddots&\vdots
\\a_{n1}&\dots&a_{nn}
\end{vmatrix}>0.
\end{equation}

It is also known that
$A$ is positive
definite if and only
if
\begin{equation}\label{ksr}
\delta_1>0,\quad
\delta_2>0,\ \dots,\
\delta_n>0,
\end{equation}
where $\delta_i$ is
the sum of all
$i\times i$ principal
minors of $A$ (a minor
is called
\emph{principal} if
its diagonal entries
lie on the diagonal of
$A$). Indeed, the
characteristic
polynomial of $A$ is
equal to
\begin{equation}\label{yys}
(x-\lambda_1)\cdots
(x-\lambda_n)=
x^n-\delta_1x^{n-1}+
\delta_2x^{n-2}-\dots
+(-1)^{n}\delta_n,
\end{equation}
and so the condition
\eqref{ksr} implies
the positivity of all
eigenvalues
$\lambda_1,\dots,
\lambda_n$ of $A$,
which ensures the
positive definiteness
of $A$ since $A$ is
unitarily similar to a
real diagonal matrix.
The conditions
\eqref{ksr} are
symmetric in the sense
of \cite{ste}:
permutations of rows
and the same
permutations of
columns of $A$ do not
change
$\delta_1,\dots,\delta_n$.

In this paper, we give
other examples of
criteria of positive
definiteness of the
form
\begin{equation}\label{kyj}
\Sigma_1>0,\quad
\Sigma_2>0,\ \dots,\
\Sigma_n=\det A>0,
\end{equation}
where each $\Sigma_i$
is a sum of some
$i\times i$ principal
minors.

In Section \ref{s1},
for each partition of
$A$ into blocks with
square diagonal
blocks, we construct a
criterion of positive
definiteness of the
form \eqref{kyj}. In
particular, the
criteria \eqref{iytd1}
and \eqref{ksr} are
obtained from the
partitions, in which
the number of diagonal
blocks is $n$ and,
respectively, $1$. We
show that the obtained
sums $\Sigma_i$ can be
used instead of the
leading principal
minors $\Delta_i$ in
other inertia
problems.

It would be
interesting to
describe all principal
minors criteria of
positive definiteness
of the form
\eqref{kyj}. In
Section \ref{s2} we
describe them for
$3\times 3$ matrices.
There are 6 such
criteria; 4 criteria
can be obtained from
the partitions into
blocks (as in Section
\ref{s1}) and the
remaining 2 criteria
are new.

This research was
inspired by Stepanov's
paper \cite{ste}, in
which the criterion of
positivity from
Theorem \ref{t1.1}(c)
was proved for real
symmetric block
matrices whose
diagonal blocks are at
most 3-by-3.

\section{Symmetric
critera of positive
definiteness for block
matrices} \label{s1}

Every $n\times n$
Hermitian complex
matrix $A=A^*$ defines
the Hermitian form
$x^*Ax$ with
$x=[x_1,\dots,x_n]^T$.
Forms $x^*Ax$ and
$x^*Bx$ are said to be
\emph{equivalent} if
their matrices $A$ and
$B$ are
*\!\emph{congruent},
i.e., $S^*AS=B$ for
some nonsingular $S$.

By Sylvester's Inertia
Law, every Hermitian
form $x^*Ax$ is
equivalent to the form
\begin{equation*}\label{lius}
\bar x_1x_1+\dots+\bar
x_px_p- \bar
x_{p+1}x_{p+1}-\dots-\bar
x_{p+q}x_{p+q},
\end{equation*}
where $p$ and $q$ do
not depend on the
method of reduction.
The numbers $p$ and
$q$ are equal to the
numbers of positive
and negative
eigenvalues of $A$
since $A$ is unitarily
*congruent to a real
diagonal matrix $D$
(i.e., $U^*AU=D$ for
some unitary $U$), see
\cite[Theorem
4.1.5]{hor}. If $r$ is
the rank of $A$ and
the leading principal
minors
$\Delta_1,\dots,\Delta_r$
are all nonzero, then
the numbers $p$ and
$q$ can be calculated
using the \emph{Jacobi
formula} \cite[Chapter
X, \S\,9, Formula
(133)]{gan}: $x^*Ax$
is equivalent to
\begin{equation}\label{liu}
\Delta_1\bar x_1x_1+
\frac{\Delta_2}{\Delta_1}\bar
x_2x_2+\dots+
\frac{\Delta_r}{\Delta_{r-1}}\bar
x_rx_r.
\end{equation}

In this section, we
consider an $n\times
n$ Hermitian matrix
$A$ partitioned into
$t\times t$ blocks
with square diagonal
blocks:
\begin{equation}\label{9aa}
A=\left[\begin{array}{c|c|c}
A_{11}&\dots&A_{1t}\\
\hline
\vdots&\ddots&\vdots \\
\hline
A_{t1}&\dots&A_{tt}
\end{array}
 \right],\qquad
\mbox{$A_{ii}$ is
$k_i$-by-$k_i$.}
\end{equation}
We say that $A$ is
\emph{block-unitarily
{\rm*}\!congruent} to
$B$ if $U^*AU=B$,
where $U$ is a direct
sum of $t$ unitary
matrices of sizes
$k_1\times
k_1,\,\dots,\,k_t\times
k_t$.

Let us denote by $A_i$
the leading principal
block submatrix of
\eqref{9aa} formed by
the first $i\times i$
blocks, i.e.,
\begin{equation}\label{iytd}
A_1=A_{11},\quad
A_2=\begin{bmatrix}
A_{11}&A_{12}\\A_{21}&A_{22}
\end{bmatrix},\ \dots,\
A_t=A.
\end{equation}
Each $(l,l)$ diagonal
entry of $A$ belongs
to some diagonal block
$A_{kk}$. Denote by
$\sigma_l$ the sum of
all $l\times l$
principal minors that
contain $A_{k-1}$ (if
$k>1$) and that are
contained in $A_{k}$.
For example, if
\[
A=\left[\begin{array}{cc|c|cc}
1&2&3&4&5\\2&3&4&5&6\\
\hline
3&4&5&6&7\\ \hline 4&5&6&7&8\\
5&6&7&8&9
\end{array}
\right],
\]
then
\begin{eqnarray*}
&\sigma_1=1+3,\qquad
\sigma_2=\Delta_2,\qquad
\sigma_3=\Delta_3,
     \\
 &\sigma_4=
\begin{vmatrix}
1&2&3&4\\2&3&4&5\\
3&4&5&6\\4&5&6&7
\end{vmatrix}
   +
\begin{vmatrix}
1&2&3&5\\2&3&4&6\\
3&4&5&7\\5&6&7&9
\end{vmatrix},
 \qquad
\sigma_5=\Delta_5=\det
A.
\end{eqnarray*}

\begin{thm} \label{t1.1}
Let $A$ be a Hermitian
matrix \eqref{9aa}
partitioned into
blocks such that the
leading principal
block submatrices
$A_1,\dots,A_{t-1}$
$($see \eqref{iytd}$)$
are nonsingular. Then

{\rm(a)} The number
$p$ of positive
eigenvalues of $A$ is
equal to the number of
sign changes
$($ignoring zeros$)$
in the sequence
\begin{equation}\label{isgt}
1,\ -\sigma_1,\
\sigma_2,\ -\sigma_3,\
\dots,\
(-1)^n\sigma_n.
\end{equation}

{\rm(b)} The number
$q$ of negative
eigenvalues of $A$ is
equal to the number of
sign changes
$($ignoring zeros$)$
in the sequence $1,\,
\sigma_1,\,
\sigma_2,\, \dots,\,
\sigma_n$.

{\rm(c)} The form
$x^*Ax$ is positive
definite if and only
if
\begin{equation}\label{kyg}
\sigma_1>0,\
\sigma_2>0,\ \dots, \
\sigma_{n}>0.
\end{equation}

{\rm(d)} The form
$x^*Ax$ is positive
semidefinite if and
only if all
$\sigma_i\ge 0$. If
all $\sigma_i\ge 0$,
then either all
$\sigma_i>0$, or all
$\sigma_i=0$, or
\begin{equation}\label{kjfq}
\sigma_1>0,\ \dots, \
\sigma_{l-1}>0,\quad
\sigma_{l}=\dots=\sigma_n=0
\end{equation}
for some $l>1$.

{\rm(e)} The number
$r:=\max\{i\,|\,\sigma_i
\ne 0\}$ is equal to
the rank of $A$. If
$\sigma_1,\dots,\sigma_{r}$
are nonzero then
$x^*Ax$ is equivalent
to
\begin{equation}\label{liuh1}
\sigma_1\bar x_1x_1+
\frac{\sigma_2}{\sigma_1}\bar
x_2x_2+\dots+
\frac{\sigma_r}{\sigma_{r-1}}\bar
x_rx_r.
\end{equation}

{\rm(f)} The numbers
$\sigma_1,\,\dots,\,\sigma_n$
are invariant with
respect to
transformations of
block-unitary
{\rm*}\!congruence
with $A$ $($in
particular, with
respect to any
permutation of rows of
$A$ within horizontal
strips and the same
permutation of its
columns$)$.
\end{thm}

\begin{proof}
We begin with a
general result on
$\sigma_i$ which will
be used in the proof
of (a)--(f). Let
$t>1$. Represent $A$
in the form
\begin{equation}\label{lkb}
A=\begin{bmatrix}
A_{t-1}&B\\B^*&A_{tt}
\end{bmatrix},\qquad
B^*=[A_{t1}|\dots|A_{t,t-1}].
\end{equation}
The size of $A_{t-1}$
is $k\times k$, where
\begin{equation}\label{po}
k:=k_1+k_2+\dots+k_{t-1}
\end{equation}
(see \eqref{9aa}). By
the assumption of the
theorem, $A_{t-1}$ is
nonsingular. Adding
linear combinations of
columns of
$A_{t-1}=A_{t-1}^*$ to
columns of $B$ and
performing the
*congruent
transformations of
rows, we reduce $A$ to
the block-diagonal
matrix
\begin{multline}\label{wevc}
\begin{bmatrix}
A_{t-1}&0\\0&A'
\end{bmatrix}:=
\begin{bmatrix}
A_{t-1}&0\\0&A_{tt}
-B^*A_{t-1}^{-1}B
\end{bmatrix}
      \\
=\begin{bmatrix}
I&0\\-(A_{t-1}^{-1}B)^*&I
\end{bmatrix}
\begin{bmatrix}
A_{t-1}&B\\B^*&A_{tt}
\end{bmatrix}
\begin{bmatrix}
I&-A_{t-1}^{-1}B\\0&I
\end{bmatrix},
\end{multline}
which is *congruent to
$A$.

If $\Delta$ is a
principal minor of $A$
and $\Delta$ contains
$A_{t-1}$, then
$\Delta$ is not
changed by
transformations
\eqref{wevc}. So
$\Delta=\Delta_k\Delta'$,
where $\Delta_k=\det
A_{t-1}$ ($k$ is
defined in \eqref{po})
and $\Delta'$ is a
principal minor of
$A'$. We have
\begin{equation}\label{ouuvg}
\sigma_{k}=\Delta_k,\
\
\sigma_{k+1}=\Delta_k\sigma'_1,
\ \ \dots, \ \
\sigma_n=
\sigma_{k+k_t}
=\Delta_k\sigma'_{k_t},
\end{equation}
where $\sigma'_j$ is
the sum of all
$j$-by-$j$ principal
minors of the matrix
$A'$.

(a) We prove the
statement (a) using
induction on $t$. Let
first $t=1$ and let
\begin{equation}\label{ldyf}
\chi_A(x)=x^n+c_1x^{n-1}+
\dots+c_n
\end{equation}
be the characteristic
polynomial of $A$.
Then
\begin{equation}\label{dou}
c_1=-\sigma_1,\
c_2=\sigma_2,\
c_3=-\sigma_3,\
\dots,\
c_n=(-1)^n\sigma_n,
\end{equation}
and the sequence
\eqref{isgt} takes the
form $1,\,
c_1,\,\dots,\,c_n$. So
the statement (a)
follows from
Descartes' Sign Rule
(\cite[\S\:55]{dic} or
\cite[Chapter 6,
\S\:4]{kos}): if all
of the roots of a
polynomial
$$f(x)=x^n+a_1x^{n-1}+
\dots+a_n\in\mathbb
R[x]$$ are real, then
the number of its
positive roots is
equal to the number of
sign changes in the
sequence of
coefficients $1,
a_1,\dots,a_n$.

Let now $t>1$. Reduce
$A$ to the form
\eqref{wevc}. By
induction hypothesis,
the statement (a)
holds for $A_{t-1}$
and $A'$. Hence, the
number $p_{t-1}$ of
positive eigenvalues
of $A_{t-1}$ is equal
to the number of sign
changes in the
sequence
$$1,\ -\sigma_1,\ \sigma_2,\
-\sigma_3,\ \dots,\
(-1)^{k}\sigma_{k},$$
and the number $p'$ of
positive eigenvalues
of $A'$ is equal to
the number of sign
changes in the
sequence
$$
1,\ -\sigma'_1,\
\sigma'_2,\
-\sigma'_3,\ \dots,\
(-1)^{k_t}\sigma'_{k_t}.
$$
In view of
\eqref{ouuvg}, the
multiplication of the
last sequence by
$(-1)^{k}\Delta_{k}$
gives the sequence
$$
(-1)^{k}\sigma_{k},\
(-1)^{k+1}\sigma_{k+1},\
\dots,\
(-1)^{n}\sigma_{n}.
$$
Therefore, the number
$p_{t-1}+p'$ of
positive eigenvalues
of (\ref{wevc}) is
equal to the number of
sign changes in the
sequence (\ref{isgt}).
This proves (a) since
by Sylvester's Inertia
Law the matrices
(\ref{wevc}) and $A$
have the same number
of positive
eigenvalues.

(b) Property (b) is
evident from property
(a) with changing $A$
by $-A$.

(c) The form $x^*Ax$
is positive definite
if and only if all the
eigenvalues of $A$ are
positive. So (c)
follows from (a).

(d) The form $x^*Ax$
is positive
semidefinite if and
only if all the
eigenvalues of $A$ are
nonnegative. So the
first statement in (d)
follows from (b).

Suppose all
$\sigma_i\ge 0$, there
exist $\sigma_i>0$,
and there exist
$\sigma_i=0$. Write
$l:=\min\{i\, |\,
\sigma_i=0\}$. Let us
prove \eqref{kjfq}
using induction on
$t$.

If $t=1$, then we
reduce $A$ by
transformations of
unitary *congruence to
a real diagonal matrix
\begin{equation}\label{lgty}
D=\diag(\lambda_1,
\dots,
\lambda_s,0,\dots,0),
\qquad \lambda_1>0,
\dots, \lambda_s>0.
\end{equation}
These transformations
do not change
$\chi_A(x)$. By
\eqref{dou}, they do
not change all
$\sigma_i$, and so we
can calculate
$\sigma_i$ using
minors of $D$ instead
of minors of $A$:

\begin{equation}\label{dyv}
\sigma_1=\sum_i
\lambda _i,\quad
\sigma_2=\sum_{i<j}
\lambda_i\lambda_j,\quad
\sigma_3=\sum_{i<j<k}
\lambda_i\lambda_j\lambda_k,\
\dots
\end{equation}
Since $\lambda_1,
\dots, \lambda_s$ are
positive, we have
\eqref{kjfq} with
$l=s$.

If $t>1$, then we
reduce $A$ to the form
\eqref{wevc}. By
induction hypothesis,
the statement (d)
holds for $A_{t-1}$
and $A'$. Since
$A_{t-1}$ is
nonsingular,
$\sigma_{k}=\Delta_k>0$,
hence all
$\sigma_{1},\dots,
\sigma_{k}$ are
positive, and so
$l>k$.

If $l=k+1$, then
$\sigma'_1=0$, and
therefore all
$\sigma'_i$ are zero.
If $l>k+1$, then
$$\sigma'_1>0,\ \dots,
\
\sigma'_{l-k-1}>0,\quad
\sigma'_{l-k}=\dots
=\sigma'_{k_t}=0.$$ In
view of \eqref{ouuvg},
this proves
\eqref{kjfq}.

(e) Let
$r:=\max\{i\,|\,\sigma_i
\ne 0\}$. Since
$A_{t-1}$ is
nonsingular,
$\sigma_k=\det
A_{t-1}\ne 0$, thus
$r\ge k$. Reduce $A$
to the form
\eqref{wevc} and
obtain \eqref{ouuvg}.
Then reduce $A'$ by
transformations of
unitary *congruence to
a real diagonal matrix
\eqref{lgty} and
obtain \eqref{dyv}
with $\sigma_i$
replaced by
$\sigma_i'$. By
\eqref{ouuvg},
\[
r=k+\max\{i\,|\,\sigma'_i
\ne 0\}=k+\rank D=
\rank A_{t-1}+\rank
A'= \rank A.
\]

If all
$\sigma_1,\dots,\sigma_{r}$
are nonzero, then the
forms $x^*Ax$ and
\eqref{liuh1} are
equivalent. Indeed,
their matrices have
the same number of
positive eigenvalues
and the same number of
negative eigenvalues
due to (a), (b), and
the equalities
$\sigma_{r+1}=\dots
=\sigma_n=0$.

(f) We use induction
on $t$. For $t=1$,
property (f) holds by
\eqref{dou} since the
coefficients of
$\chi_A(x)$ are
invariant with respect
to similarity
transformations with
$A$. For $t>1$,
consider $\widetilde
A:=U^*AU$, where
$U=U_1\oplus\dots\oplus
U_t$ and each $U_i$ is
a $k_i\times k_i$
unitary matrix. The
sums $\sigma_i$ were
defined for $A$;
denote by
$\widetilde\sigma_1,\dots,
\widetilde\sigma_n$
the corresponding sums
for $\widetilde A$.
Partition $A$ into
blocks as in
(\ref{lkb}) and
partition $\widetilde
A$ analogously:
\begin{equation*}\label{lkbs}
\widetilde
A=\begin{bmatrix}
\widetilde
A_{t-1}&\widetilde
B\\\widetilde
B^*&\widetilde A_{tt}
\end{bmatrix},\qquad
\widetilde
B^*:=[\widetilde
A_{t1}|\dots|\widetilde
A_{t,t-1}].
\end{equation*}

Let $V:=U_1\oplus
\dots\oplus U_{t-1}$.
Then $\widetilde
A_{t-1}= V^*A_{t-1}V$.
By induction
hypothesis, property
(f) holds for
$A_{t-1}$, that is, $
 \sigma_1=\widetilde \sigma_1,\
 \dots,\ \sigma_k=\widetilde
 \sigma_k,
$ where $k$ was
defined in \eqref{po}.
It remains to prove
that
\begin{equation}\label{ljd}
 \sigma_{k+1}=\widetilde \sigma_{k+1},\
 \dots,\ \sigma_n=\widetilde
 \sigma_n.
\end{equation}
Since $U=(V\oplus
I_{k_t})(I_k\oplus
U_t)$, the
transformation
$A\mapsto U^*AU$ is
the composition of two
transformations:
$A\mapsto (V\oplus
I_{k_t})^*A(V\oplus
I_{k_t})$ and
$A\mapsto (I_k\oplus
U_t)^*A(I_k\oplus
U_t)$. The first
transformation does
not change
$\sigma_{k+1},
\dots,\sigma_n$ since
it does not change
every minor of $A$
containing $A_{t-1}$.
It remains to prove
\eqref{ljd} for the
second transformation.

Thus we can suppose
that $V=I_k$. Then
\[
\widetilde A=U^*AU=
\begin{bmatrix}
A_{t-1}&BU_t
\\U_t^*B^*&U_t^*A_{tt}U_t
\end{bmatrix}.
\]
Reduce $A$ to the form
(\ref{wevc}) and
$\widetilde A$ to the
form
\begin{multline}\label{wcevc}
\begin{bmatrix}
A_{t-1}&0\\0&\widetilde
A'
\end{bmatrix}:=
\begin{bmatrix}
A_{t-1}&0\\0&U_t^*A'U_t
\end{bmatrix}
      \\ 
=\begin{bmatrix}
I&0\\-(A_{t-1}^{-1}BU_t)^*&I
\end{bmatrix}
\begin{bmatrix}
A_{t-1}&BU_t
\\U_t^*B^*&U_t^*A_{tt}U_t
\end{bmatrix}
\begin{bmatrix}
I&-A_{t-1}^{-1}BU_t\\0&I
\end{bmatrix},
\end{multline}
where $A'=A_{tt}
-B^*A_{t-1}^{-1}B$ was
defined in
(\ref{wevc}).

Each $\sigma_{i}$ with
$i\ge k$ has the form
\eqref{ouuvg}.
Analogously, each
$\widetilde\sigma_{i}$
with $i\ge k$ has the
form
\begin{equation*}
\widetilde
\sigma_{k}=\Delta_k,\
\
\widetilde\sigma_{k+1}=\Delta_k
\widetilde\sigma'_1, \
\ \dots, \ \
\widetilde\sigma_{k+r-1}=\Delta_k
\widetilde\sigma'_{r-1},
\end{equation*}
where
$\widetilde\sigma'_j$
is the sum of all
$j$-by-$j$ principal
minors of the matrix
$\widetilde A'$. Since
the matrices $A'$ and
$\widetilde
A'=U_t^*A'U_t$ are
similar, $\sigma'_j
=\widetilde\sigma'_j$
for all $j$. This
proves \eqref{ljd}.
\end{proof}

\section{Principal minors
criteria for $3\times
3$ matrices}
\label{s2}

For each $n\times n$
matrix, denote by
\[
P_{i_1i_2\ldots
i_k},\qquad 1\le
i_1<i_2<\dots< i_k\le
n,
\]
its $k\times k$
principal minor that
lies on the
intersection of rows
$i_1,i_2,\dots, i_k$
with columns
$i_1,i_2,\dots, i_k$.
Let \eqref{kyj} be a
system of
inequalities, in which
every $\Sigma_k$ is a
sum of some
$P_{i_1i_2\ldots i_k}$
with distinct
$(i_1,i_2,\ldots,
i_k)$. We say that
\eqref{kyj}
\emph{ensures positive
definiteness} if every
$n\times n$ Hermitian
matrix is positive
definite if and only
if it satisfies
\eqref{kyj}.

\begin{thm} \label{t33}
{\rm(a)} Each of the
following systems of
inequalities ensures
positive definiteness
of $3$-by-$3$
Hermitian matrices:
\begin{itemize}
  \item[\rm(i)]
$P_1>0$, \ $P_{12}>0$,
\ $P_{123}>0$;

  \item[\rm(ii)]
$P_1>0$, \ $P_{12}
+P_{13}>0$, \
$P_{123}>0$;

  \item[\rm(iii)]
$P_1+P_2>0$, \
$P_{12}>0$, \
$P_{123}>0$;

  \item[\rm(iv)]
$P_1+P_2>0$, \
$P_{12}+ P_{13}+
P_{23}>0$, \
$P_{123}>0$;

  \item[\rm(v)]
$P_1+P_2 +P_3>0$, \
$P_{12}+ P_{13}>0$, \
$P_{123}>0$;

  \item[\rm(vi)]
$P_1+P_2 +P_3>0$, \
$P_{12}+ P_{13}+
P_{23}>0$, \
$P_{123}>0$.
\end{itemize}
Systems {\rm(i)},
{\rm(ii)}, {\rm(iii)},
and {\rm(vi)} have the
form \eqref{kyg} with
respect to the
partitions
\begin{equation}\label{iyu}
\left[\begin{array}{c|c|c}
\ci&\ci&\ci\\\hline
\ci&\ci&\ci\\\hline
\ci&\ci&\ci
\end{array}
\right],
   \quad
\left[\begin{array}{c|cc}
\ci&\ci&\ci\\\hline
\ci&\ci&\ci\\
\ci&\ci&\ci
\end{array}
\right],
      \quad
\left[\begin{array}{cc|c}
\ci&\ci&\ci\\
\ci&\ci&\ci\\\hline
\ci&\ci&\ci
\end{array}
\right],\quad
\left[\begin{array}{ccc}
\ci&\ci&\ci\\
\ci&\ci&\ci\\
\ci&\ci&\ci
\end{array}
\right].
\end{equation}

{\rm(b)} If a system
\eqref{kyj} with $n=3$
ensures positive
definiteness, then it
can be obtained from
one of {\rm(i)--(vi)}
by a permutation of
the indexing set
$\{1,2,3\}$.
\end{thm}

\begin{proof}
If a $3$-by-$3$
Hermitian matrix is
positive definite,
then its principal
minors are positive,
and so it satisfies
each of systems
(i)--(vi).

Let
\begin{equation}\label{jdo}
\Sigma_1>0,\quad
\Sigma_2>0,\quad
P_{123}>0
\end{equation}
be a system of the
form \eqref{kyj} with
$n=3$. For each
substitution $\sigma$
on the indexing set
$\{1,2,3\}$, we define
the system
\begin{equation}\label{jdo1}
\Sigma_1^{\sigma}>0,\quad
\Sigma_2^{\sigma}>0,\quad
P_{123}>0
\end{equation}
obtained from
\eqref{jdo} by
replacement of all the
summands $P_i$ and
$P_{ij}$ of $\Sigma_1$
and $\Sigma_2$ with
$P_{\sigma (i)}$ and
$P_{\sigma (i)\sigma
(j)}$. A $3\times 3$
Hermitian matrix $A$
satisfies \eqref{jdo}
if and only if the
matrix $A^{\sigma}$
obtained by the
corresponding
permutations of rows
and columns satisfies
\eqref{jdo1}. Hence
\eqref{jdo} ensures
positive definiteness
if and only if the
same holds for
\eqref{jdo1}.

Each system of the
form \eqref{jdo}
determined up to
substitutions $\sigma$
is presented by one of
the rows of the
following table:
\[
\begin{array}{|l|l|c|}
\hline
 P_1&
P_{12} & \text{(i)}\\
& P_{23}
& \diag(1,-1,-1)\\
& P_{12}+P_{13}
& \text{(ii)}\\
& P_{12}+P_{23}
& \diag(1,-1,-2)\\
& P_{12}+P_{13}
+P_{23}
& \diag(1,-2,-3)\\
        \hline
 P_1+P_2&
P_{12} & \text{(iii)}\\
& P_{13}
& \diag(-1,2,-1)\\
& P_{12}+P_{13}
& \diag(1,-2,-1)\\
& P_{12}+P_{13}
+P_{23}
& \text{(iv)}\\
\hline
 P_1+P_2
 +P_3&
P_{12} &
\diag(-1,-1,3)\\
& P_{12}+ P_{13}
& \text{(v)}\\
& P_{12}+P_{13}
+P_{23}
& \text{(vi)}\\
\hline
\end{array}
\]
The first two entries
of the row are
$\Sigma_1$ and
$\Sigma_2$, and the
last entry is either a
matrix that is not
positive definite but
fulfils $\Sigma_1>0$,
$\Sigma_2>0$,
$P_{123}>0$ (which
means that the system
does not ensure
positive definiteness)
or the number
((i)--(vi)) of the
corresponding system
in Theorem \ref{t33}.

It remains to prove
that each of the
systems (i)--(vi)
ensures positive
definiteness. This is
true for (i), (ii),
(iii), and (vi) due to
Theorem \ref{t1.1}(c)
applied to $3\times 3$
matrices partitioned
as in \eqref{iyu}.

Let a $3\times 3$
Hermitian matrix $A$
satisfy condition
(iv). By a suitable
transformation
$(U\oplus I_1)^{-1}A
(U\oplus I_1)$ with
unitary $U$, we reduce
$A$ to the form
\begin{equation}\label{ljv}
\begin{bmatrix}
a&0&\bar x\\
0&b&\bar y\\
x&y&c
\end{bmatrix}.
\end{equation}
This similarity
transformation does
not change the
left-hand sides of the
inequalities (iv)
since $P_1+P_2$ is the
trace of the leading
principal $2\times 2$
submatrix (whose
determinant is
$P_{12}$), and
$P_{12}+ P_{13}+
P_{23}$ is a
coefficient of the
characteristic
polynomial of $A$ (see
\eqref{yys}).
Therefore, the matrix
\eqref{ljv} fulfils
(iv):
\begin{equation}\label{ljfv}
a+b>0,\quad
ab+(ac-|x|^2)
+(bc-|y|^2)>0,\quad
abc-|x|^2b-|y|^2a>0.
\end{equation}
If $c<0$ then $ab>0$
by the first and the
second inequalities in
\eqref{ljfv}; since
$a+b>0$, we have $a>0$
and $b>0$, which
contradicts the third
inequality in
\eqref{ljfv}. Thus
$c\ge 0$, $a+b+c>0$,
$A$ satisfies (vi),
and so it is positive
definite.

Let a $3\times 3$
Hermitian matrix $A$
satisfy condition (v).
By a suitable
transformation
$(I_1\oplus U)^{-1}A
(I_1\oplus U)$ with
unitary $U$, we reduce
$A$ to the form
\begin{equation}\label{ljf}
\begin{bmatrix}
a&\bar x&\bar y\\
x&b&0\\
y&0&c
\end{bmatrix}.
\end{equation}
This similarity
transformation does
not change $P_{23}$
and $P_{12}+ P_{13}+
P_{23}$, hence it
preserves $P_{12}+
P_{13}$. Therefore,
the matrix \eqref{ljf}
fulfils (v):
\begin{equation*}\label{ljs}
a+b+c>0,\quad
ab-|x|^2+ac-|y|^2>0,\quad
abc-|x|^2c-|y|^2b>0.
\end{equation*}
Since
$a(b+c)>|x|^2+|y|^2$,
$a\ne 0$. If $a<0$
then $b+c<0$, which
contradicts $a+b+c>0$.
Thus $a>0$, $A$
satisfies (ii), and so
it is positive
definite.
\end{proof}
\medskip

\subsection*{Acknowledgment}
The authors are
greatly indebted to
Professor Roger Horn
for many helpful
comments. The authors
also wish to express
their gratitude to the
referee for suggesting
new problems, which
were partially solved
in the revised
version.

\end{document}